\documentclass[12pt, a4paper]{amsart}

\usepackage{enumerate}
\usepackage{tikz}
\usetikzlibrary{fit,topaths,calc}

\newtheorem{firstthm}{Proposition}

\newtheorem{lemma}[firstthm]{Lemma}
\newtheorem{prop}[firstthm]{Proposition}
\newtheorem{conj}[firstthm]{Conjecture}
\newtheorem{construct}[firstthm]{Construction}

\def\COMMENT#1{}
\let\COMMENT=\footnote

\date{\today}
\title{Two constructions relating to conjectures of Beck on positional games}
\author{Fiachra Knox}

\begin{document}

\vspace*{-0.8cm}
\begin{abstract}
In this paper, we construct two hypergraphs which exhibit the following properties.
We first construct a hypergraph $G_{CP}$ and show that Breaker wins the Maker-Breaker game on $G_{CP}$, but Chooser wins the Chooser-Picker game on $G_{CP}$.
This disproves an (informally stated) conjecture of Beck.
Our second construction relates to Beck's Neighbourhood Conjecture, which (in its weakest form) states that there exists $c > 1$ such that Breaker wins the Maker-Breaker game on any $n$-uniform hypergraph $G$ of maximum degree at most $c^n$.
We consider the case $n=4$ and construct a $4$-graph $G_4$ with maximum vertex degree $3$, such that Maker wins the Maker-Breaker game on $G_4$.
This answers a question of Leader.
\end{abstract}
\maketitle
\vspace*{-0.6cm}

\section{Introduction}

Positional games have been the object of systematic study since the 1960's.
One of the most well-known and well-studied of these games is the \emph{Maker-Breaker game}.
Given a (finite) hypergraph $G$, the Maker-Breaker game on $G$ is a $2$-player game which is defined as follows. 
Two players, Maker and Breaker, each take turns to claim a (previously unclaimed) vertex of $G$. 
Maker's goal is to claim all of the vertices of some edge of $G$; if he accomplishes this at any point, then he wins the game.
However, if all of the verties of $G$ have been claimed and Maker has not achieved his goal, then Breaker wins.
An equivalent formulation of Breaker's goal is that he wins if he claims at least one vertex of every edge of $G$.

Major early breakthroughs in the area of positional games came in the form of the Hales-Jewett Theorem~\cite{HJ63} and the Erd\H{o}s-Selfridge Theorem~\cite{ES73}, the latter of which states that the condition $\sum_{A \in G} 2^{-|A|} < 1/2$ is a sufficient condition for Breaker to win.
A thorough study of these games has been made by Beck over a lengthy period, and~\cite{BeckBook} provides a comprehensive introduction to positional games and their applications.
In general positional games are perfect-information games and involve no randomness.
Thus in theory, perfect play is possible and the outcome of the game depends only upon the kind of game, the hypergraph upon which the game is played and who (if applicable) has the first move.
(Unless stated otherwise, we assume in what follows that Maker goes first in the Maker-Breaker game.)
When a player, e.g. Maker, has a winning strategy for a certain game then we say that the game is a \emph{win for Maker}, or simply that Maker `wins' the game.

\subsection{Chooser-Picker games}

In addition to the Maker-Breaker game, many other positional games have been studied.
One such game is the \emph{Chooser-Picker game}, which is played between two players, Chooser and Picker.
As in the Maker-Breaker game, the first player (Chooser) wins if he claims every vertex of some edge and the second player (Picker) wins otherwise.
However, instead of taking turns at claiming vertices as in the Maker-Breaker game, the vertices are allocated by repeating the following procedure: On each turn, Picker selects a pair $\{x, y\}$ of unclaimed vertices to offer to Chooser. 
Chooser then claims one of the offered vertices; however, if he claims $x$ then Picker is allowed to claim $y$, and vice-versa.
If a position is reached in which only one unclaimed vertex remains, then Chooser automatically claims the final vertex.

A natural question is whether there is any relationship between the Maker-Breaker and Chooser-Picker games when the two are played on the same hypergraph $G$.
For some hypergraphs, Picker has a much easier time playing the Chooser-Picker game than Breaker has playing the Maker-Breaker game on the same hypergraph.
For example, consider a $k$-uniform hypergraph (hereafter referred to as a \emph{$k$-graph}) $K^k_{n, \ldots, n}$ which is complete $k$-partite on vertex classes of size $n \geq 2$.
Maker has a simple winning strategy for the Maker-Breaker game on $K^k_{n, \ldots, n}$: Whenever Breaker first claims a vertex of a vertex class $V$, Maker responds by claiming any unclaimed vertex of $V$ (if he has not already done so).
However, if $k \geq 2^{n}$ then Picker has a winning strategy for the Chooser-Picker game on $K^k_{n, \ldots, n}$.
Indeed, Picker can guarantee claiming every vertex of some vertex class $V$, which suffices for a win since every edge of $K^k_{n, \ldots, n}$ contains some vertex of $V$.

Hypergraphs on which Picker has a \emph{harder} time than Breaker seem to be far less common.
The general trend of Chooser-Picker games being more favourable to Picker (and less favourable to Chooser) than Maker-Breaker games was noted by Beck, and the following precise conjecture was formulated in~\cite{CMP09}.

\begin{conj}\cite{CMP09} \label{CPConj} Let $G$ be a hypergraph and suppose that Breaker wins the Maker-Breaker game on $G$, where Maker goes first. Then Picker wins the Chooser-Picker game on $G$.
\end{conj}

In this paper we give a counterexample to Conjecture~\ref{CPConj}; that is, we construct a hypergraph $G_{CP}$ such that Chooser wins the Chooser-Picker game on $G_{CP}$, but Breaker wins the Maker-Breaker game on the same graph.
A natural question to ask is what happens in the Picker-Chooser game; that is, when Picker's goal is to claim every vertex of some edge of $G$ and Chooser's goal is to prevent this.
A similar conjecture to Conjecture~\ref{CPConj} concerning the Picker-Chooser game was formulated in~\cite{CMP09} and shown to be equivalent to Conjecture~\ref{CPConj}.

Even though Conjecture~\ref{CPConj} turns out to be false, it is known that in some special cases the conjecture does hold.
Csernenszky, M\'andity and Pluh\'ar~\cite{CMP09} showed that a stronger version of the Erd\H{o}s-Selfridge condition on $G$ (with a smaller quantity in place of $1/2$) implies that Picker wins the Chooser-Picker game on $G$.
They conjectured that the Erd\H{o}s-Selfridge condition itself suffices; that is, if $\sum_{e \in G} 2^{-|e|} < 1/2$ then Picker wins the Chooser-Picker game on $G$.
Very recently, this was proved by Bednarska-Bzd\c{e}ga~\cite{MB-Bprep}.

\subsection{The Neighbourhood conjecture} \label{G4Intro}

The Neighbourhood conjecture, proposed by Beck~\cite[Open Problem 9.1]{BeckBook} in several forms, is a central problem in the theory of positional games.
In its strongest form, the Neighbourhood conjecture states that if $G$ is an $n$-graph and every edge of $G$ intersects fewer than $2^{n-1}$ other edges, then Breaker wins the Maker-Breaker game on $G$.
This version of the conjecture was disproved by Gebauer~\cite{HG}.
However perhaps the most interesting version is the following weaker statement: There exists $c > 1$ such that Breaker wins the Maker-Breaker game on any $n$-graph with maximum (vertex) degree at most $c^n$.

For a given $n$, let $f(n)$ be the smallest integer for which there exists an $n$-graph $G$ of maximum degree at most $f(n)$ such that Maker wins the Maker-Breaker game on $G$. 
A well-known pairing argument due to Hales and Jewett~\cite{HJ63} shows that $f(n) > n/2$ for each $n$.
Surprisingly, this remains the best known lower bound for $f(n)$.
On the other hand the best known upper bound on $f(n)$, which follows from a result of Gebauer, Szab\'o and Tardos~\cite{GST11}, is within a constant factor of $2^n/n$.
Given the apparent intractability of determining the general behaviour of $f(n)$, we examine the values of $f(n)$ for small $n$.

\begin{figure}
\begin{center}
\begin{tikzpicture}[label distance=2]

\node[label=above:$v_1$] (v1) at (3.5, 3.5) {};
\node[label=left:$v_2$] (v2) at (1.5, 1.5) {};
\node[label=right:$v_3$] (v3) at (5.5, 1.5) {};
\node (v4) at (2, 1) {};
\node (v5) at (2, 2) {};
\node (v6) at (2, 3) {};
\node (v7) at (3, 1) {};
\node (v8) at (3, 2) {};
\node (v9) at (3, 3) {};
\node (v10) at (4, 1) {};
\node (v11) at (4, 2) {};
\node (v12) at (4, 3) {};
\node (v13) at (5, 1) {};
\node (v14) at (5, 2) {};
\node (v15) at (5, 3) {};

\foreach \i in {1,2,...,15} {
        \fill (v\i) circle (0.1);
}

\draw[rounded corners=8] ($(v4) + (-0.2828, -0.2828)$) -- ($(v6) + (-0.2828, 0.2828)$) -- ($(v6) + (0.2828, 0.2828)$) -- ($(v4) + (0.2828, -0.2828)$) -- cycle;
\draw[rounded corners=8] ($(v7) + (-0.2828, -0.2828)$) -- ($(v9) + (-0.2828, 0.2828)$) -- ($(v9) + (0.2828, 0.2828)$) -- ($(v7) + (0.2828, -0.2828)$) -- cycle;
\draw[rounded corners=8] ($(v10) + (-0.2828, -0.2828)$) -- ($(v12) + (-0.2828, 0.2828)$) -- ($(v12) + (0.2828, 0.2828)$) -- ($(v10) + (0.2828, -0.2828)$) -- cycle;
\draw[rounded corners=8] ($(v13) + (-0.2828, -0.2828)$) -- ($(v15) + (-0.2828, 0.2828)$) -- ($(v15) + (0.2828, 0.2828)$) -- ($(v13) + (0.2828, -0.2828)$) -- cycle;

%\draw[rounded corners=8] ($(v1) + (-0.2828, 0.2828)$) -- ($(v15) + (0.2828, 0.2828)$) -- ($(v15) + (0.2828, -0.2828)$) -- ($(v1) + (-0.2828, -0.2828)$) -- cycle;
%\draw[rounded corners=8] ($(v1) + (0.2828, 0.2828)$) -- ($(v6) + (-0.2828, 0.2828)$) -- ($(v6) + (-0.2828, -0.2828)$) -- ($(v1) + (0.2828, -0.2828)$) -- cycle;

\draw[rounded corners=8] ($(v1) + (-0.4,0)$) -- ($(v12) + (-0.1112,-0.2828)$) -- ($(v15) + (0.2828,-0.2828)$) -- ($(v15) + (0.2828,0.2828)$) -- ($(v12) + (0.1112,0.2828)$) -- ($(v1) + (0,0.4)$) -- cycle;
\draw[rounded corners=8] ($(v1) + (0.4,0)$) -- ($(v9) + (0.1112,-0.2828)$) -- ($(v6) + (-0.2828,-0.2828)$) -- ($(v6) + (-0.2828,0.2828)$) -- ($(v9) + (-0.1112,0.2828)$) -- ($(v1) + (0,0.4)$) -- cycle;

\draw[rounded corners=8] ($(v2) + (-0.4,0)$) -- ($(v4) + (-0.1112,-0.2828)$) -- ($(v7) + (0.2828,-0.2828)$) -- ($(v7) + (0.2828,0.2828)$) -- ($(v4) + (0.1112,0.2828)$) -- ($(v2) + (0,0.4)$) -- cycle;
\draw[rounded corners=8] ($(v2) + (-0.4,0)$) -- ($(v5) + (-0.1112,0.2828)$) -- ($(v8) + (0.2828,0.2828)$) -- ($(v8) + (0.2828,-0.2828)$) -- ($(v5) + (0.1112,-0.2828)$) -- ($(v2) + (0,-0.4)$) -- cycle;

\draw[rounded corners=8] ($(v3) + (0.4,0)$) -- ($(v13) + (0.1112,-0.2828)$) -- ($(v10) + (-0.2828,-0.2828)$) -- ($(v10) + (-0.2828,0.2828)$) -- ($(v13) + (-0.1112,0.2828)$) -- ($(v3) + (0,0.4)$) -- cycle;
\draw[rounded corners=8] ($(v3) + (0.4,0)$) -- ($(v14) + (0.1112,0.2828)$) -- ($(v11) + (-0.2828,0.2828)$) -- ($(v11) + (-0.2828,-0.2828)$) -- ($(v14) + (-0.1112,-0.2828)$) -- ($(v3) + (0,-0.4)$) -- cycle;

\end{tikzpicture}
\caption{The graph $G_3$.} \label{G3figure}
\end{center}
\end{figure}
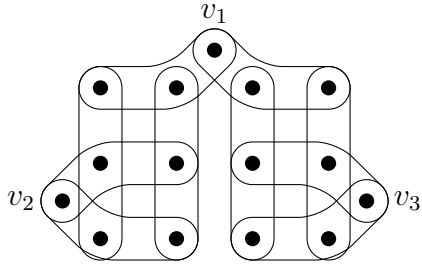

It is easy to see that $f(2) = 2$. 
The graph $G_3$ (see Figure~\ref{G3figure}), a $3$-graph of maximum degree $2$, is a win for Maker; this demonstrates that $f(3) = 2$ also. 
Maker's winning strategy is to claim $v_1$, and then if Breaker plays to the right of $v_1$ he claims $v_2$; otherwise, he claims $v_3$. In either case he wins quickly regardless of Breaker's next move.

We can derive from $G_3$ a $4$-graph which has maximum degree $4$, and on which Maker wins, as follows: For each edge $e \in G_3$, add two vertices $x_e$ and $y_e$ and replace $e$ by the two edges $e \cup \{x_e\}$ and $e \cup \{y_e\}$. 
To win, Maker first claims an edge of $G_3$ and then claims one of $x_e$ and $y_e$, making appropriate modifications if Breaker claims some $x_e$ or $y_e$ prematurely.
We therefore deduce that $f(4) \leq 4$.
On the other hand $f(4) > 4/2 = 2$ by the pairing argument, leaving only $3$ and $4$ as possible values for $f(4)$.
A question posed by Imre Leader at the Workshop on Probabilistic techniques in Graph Theory, University of Birmingham (March 25, 2012) is whether there exists a $4$-graph with maximum degree $3$ on which Maker wins the Maker-Breaker game. 

In this paper we exhibit a $4$-graph $G_4$ which answers this question in the affirmative, proving that $f(4) = 3$.
Our counterexample is quite large, with $472$ vertices and $331$ edges. 
The exact value for $f(n)$ is unknown for any $n \geq 5$.

\subsection{Notation}

We denote by $[r]$ the set $\{1, 2, \ldots, r\}$. 
For simplicity we refer to the maximum vertex degree of a graph $G$ simply as its maximum degree (other kinds of degree will not be considered in this paper).
Finally, we say that a player has \emph{claimed} an edge $e$ when that player has claimed every vertex of $e$.

\subsection{Acknowledgements} I would like to thank Malgorzata Bednarska-Bzd\c{e}ga for introducing the `asymmetry' condition described at the start of Section~\ref{CPSection}, which any counterexample to Conjecture~\ref{CPConj} must satisfy; and also for many helpful discussions.
I would also like to thank Deryk Osthus and Daniela K\"{u}hn for their comments on the Sections~\ref{G4Intro} and~\ref{G4Section}.

\section{Preliminaries}

Before presenting the graphs $G_{CP}$ and $G_4$, we establish some preliminaries and give the previously mentioned pairing argument of Hales and Jewett, which shows that $f(n) > n/2$.
Given an $n$-graph $G$ with maximum degree at most $n/2$, form an auxiliary bipartite $2$-graph $B$ on vertex classes $V_1 = V(G)$ and $V_2$, where $V_2$ is composed of \emph{two} copies of $E(G)$. 
We join $v \in V_1$ to $e \in V_2$ in $B$ if and only if $v \in e$ in $G$.
Since $G$ is $n$-uniform, the degree of every vertex in $V_2$ is $n$.
Further, the degree of a vertex in $V_1$ is twice that of the corresponding vertex in $G$, and hence is at most $n$.
It is not hard to see that $B$ satisfies Hall's condition and so Hall's theorem implies that $B$ contains a matching which covers every vertex of $V_2$.
For each $e \in E(G)$, let $x_e$ and $y_e$ be the vertices which were matched to copies of $e$ in $V_2$.
Now Breaker's winning strategy is, whenever Maker claims one of $x_e$ and $y_e$, to immediately claim the other.
This ensures that Maker cannot claim any edge $e$, since he would have to claim both $x_e$ and $y_e$.

We now establish some simple lemmas which will prove useful in analysing the Maker-Breaker and Chooser-Picker games on $G_{CP}$.
We begin by considering the Maker-Breaker game.
Intuitively, the presence of vertices of low degree in an edge makes it easier for Breaker to block that edge.
This intuition is captured in the following two lemmas.

\begin{lemma} \label{Reduction1} Let $G$ be a hypergraph and let $e$ be an edge of $G$ which contains two vertices $x$ and $y$, each of degree $1$. Suppose that Breaker wins the Maker-Breaker game on $G - \{x, y\}$, where Maker goes first. Then Breaker also wins the Maker-Breaker game on $G$, where Maker goes first.
\end{lemma}

\proof Breaker has the following winning strategy: If Maker claims one of $\{x, y\}$ then Breaker responds by claiming the other; this ensures that Maker can never claim $e$. Otherwise, Breaker follows his winning strategy for $G - \{x, y\}$; this ensures that Maker can never claim any other edge of $G$.  \endproof

\begin{lemma} \label{Reduction2} Let $G$ be a hypergraph and let $e = \{x, y\}$ be an edge of $G$ of size $2$. Suppose that $x$ has degree $1$ and that Maker wins the Maker-Breaker game on $G$. Then Maker has a winning strategy which begins by claiming $y$. \end{lemma} 

\proof Consider any winning strategy for Maker. Then we construct a new strategy as follows: Maker first claims $y$, and if Breaker does not respond by claiming $x$ then Maker wins immediately by doing so. 
If Breaker does respond at $x$, Maker simply follows his original strategy. If at any point Maker's strategy calls for him to claim one of $x$ and $y$, Maker plays as if Breaker had responded at the other. 

At the end of the game, the original strategy would have allowed Maker to claim an edge $e'$ of $G$. But $e' \neq e$, since by pretending that Breaker claimed either $x$ or $y$ we ensured that the original strategy could not allow Maker to claim $e$ given Breaker's play. So $e'$ cannot contain $x$. 
Further, if $e'$ contains $y$ then Maker still wins using the new strategy since he really did claim $y$. 
If $e'$ does not contain $x$ or $y$ then Maker will claim $e'$, since he only deviates from the original winning strategy when it calls for him to claim $x$ or $y$. \endproof

When analysing the Chooser-Picker game, we will require the following result from~\cite{CMP09}.
The central idea of the proof is one that we will frequently use later on: In the course of the Chooser-Picker game, Chooser can never afford to allow Picker to claim all but one vertex of any edge unless he claims the last vertex on the same turn.
If he does, then he will eventually be forced to offer the last vertex and so from then on he can only delay Picker's win.

\begin{lemma}\cite[Lemma 9]{CMP09} \label{Reduction3} In the Chooser-Picker game, if there is any winning set with no vertices claimed by Picker and exactly two unclaimed elements $x$ and $y$, and Picker has a winning strategy, then he has a winning strategy starting by offering $\{x, y\}$. \end{lemma}

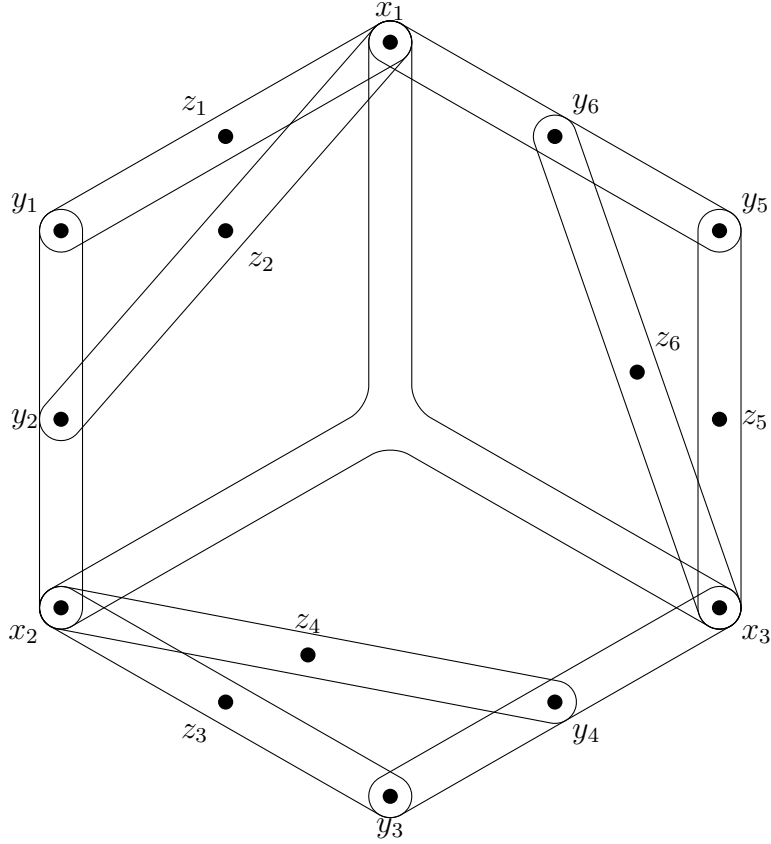
\begin{figure}
\begin{center}
\begin{tikzpicture}

\node[label=above:$x_1$] (x1) at (90: 5) {};
\node[label=150:$y_1$] (y1) at (150: 5) {};
\node[label=210:$x_2$] (x2) at (210: 5) {};
\node[label=270:$y_3$] (y3) at (270: 5) {};
\node[label=330:$x_3$] (x3) at (330: 5) {};
\node[label=30:$y_5$] (y5) at (30: 5) {};

\node[label=120:$z_1$] (z1) at (barycentric cs:x1=0.5,y1=0.5) {};
\node[label=180:$y_2$] (y2) at (barycentric cs:x2=0.5,y1=0.5) {};
\node[label=240:$z_3$] (z3) at (barycentric cs:x2=0.5,y3=0.5) {};
\node[label=300:$y_4$] (y4) at (barycentric cs:x3=0.5,y3=0.5) {};
\node[label=0:$z_5$] (z5) at (barycentric cs:x3=0.5,y5=0.5) {};
\node[label=60:$y_6$] (y6) at (barycentric cs:x1=0.5,y5=0.5) {};

\node[label=below right:$z_2$] (z2) at (barycentric cs:x1=0.5,y2=0.5) {};
\node[label=above:$z_4$] (z4) at (barycentric cs:x2=0.5,y4=0.5) {};
\node[label=60:$z_6$] (z6) at (barycentric cs:x3=0.5,y6=0.5) {};

\foreach \i in {1,2,3} {
        \fill (x\i) circle (0.1);
}

\foreach \i in {1,2,...,6} {
        \fill (y\i) circle (0.1);
}

\foreach \i in {1,2,...,6} {
        \fill (z\i) circle (0.1);
}

\draw[rounded corners=8] ($(x2) + (-0.2828, -0.2828)$) -- ($(y1) + (-0.2828, 0.2828)$) -- ($(y1) + (0.2828, 0.2828)$) -- ($(x2) + (0.2828, -0.2828)$) -- cycle;
\draw[rounded corners=8] ($(x3) + (-0.2828, -0.2828)$) -- ($(y5) + (-0.2828, 0.2828)$) -- ($(y5) + (0.2828, 0.2828)$) -- ($(x3) + (0.2828, -0.2828)$) -- cycle;
\draw[rounded corners=8] ($(x1) + (0.3864, -0.1035)$) -- ($(y1) + (-0.1035, -0.3864)$) -- ($(y1) + (-0.3864, 0.1035)$) -- ($(x1) + (0.1035, 0.3864)$) -- cycle;
\draw[rounded corners=8] ($(x3) + (0.3864, -0.1035)$) -- ($(y3) + (-0.1035, -0.3864)$) -- ($(y3) + (-0.3864, 0.1035)$) -- ($(x3) + (0.1035, 0.3864)$) -- cycle;
\draw[rounded corners=8] ($(x2) + (-0.3864, -0.1035)$) -- ($(y3) + (0.1035, -0.3864)$) -- ($(y3) + (0.3864, 0.1035)$) -- ($(x2) + (-0.1035, 0.3864)$) -- cycle;
\draw[rounded corners=8] ($(x1) + (-0.3864, -0.1035)$) -- ($(y5) + (0.1035, -0.3864)$) -- ($(y5) + (0.3864, 0.1035)$) -- ($(x1) + (-0.1035, 0.3864)$) -- cycle;

\draw[rounded corners=8] ($(x1) + (0.3990, 0.0286)$) -- ($(y2) + (0.0286, -0.3990)$) -- ($(y2) + (-0.3990, -0.0286)$) -- ($(x1) + (-0.0286, 0.3990)$) -- cycle;
\draw[rounded corners=8, transform canvas={rotate around={120:(0,0)}}] ($(x1) + (0.3990, 0.0286)$) -- ($(y2) + (0.0286, -0.3990)$) -- ($(y2) + (-0.3990, -0.0286)$) -- ($(x1) + (-0.0286, 0.3990)$) -- cycle;
\draw[rounded corners=8, transform canvas={rotate around={240:(0,0)}}] ($(x1) + (0.3990, 0.0286)$) -- ($(y2) + (0.0286, -0.3990)$) -- ($(y2) + (-0.3990, -0.0286)$) -- ($(x1) + (-0.0286, 0.3990)$) -- cycle;

\draw[rounded corners=8] ($(x1) + (-0.2828, 0.2828)$) -- (-0.2828, 0.1633) -- ($(x2) + (-0.3864, 0.1035)$) -- ($(x2) + (-0.1035, -0.3864)$) -- (0, -0.3266) -- ($(x3) + (0.1035, -0.3864)$) -- ($(x3) + (0.3864, 0.1035)$) -- (0.2828, 0.1633) -- ($(x1) + (0.2828, 0.2828)$) -- cycle;

\end{tikzpicture}
\caption{The graph $G_{CP}$ from Construction~\ref{CPConstruction}.} \label{CPfigure}
\end{center}
\end{figure}

\section{A counterexample to Conjecture~\ref{CPConj}} \label{CPSection}

Our counterexample is constructed as follows (see Figure~\ref{CPfigure}).

\begin{construct} \label{CPConstruction}
Let $G_{CP}$ be a $3$-graph on the vertex set $X \cup Y \cup Z$, where $X = \{x_1, x_2, x_3\}$, $Y = \{y_1, \ldots, y_6\}$ and $z = \{z_1, \ldots, z_6\}$, whose edge set consists of
\begin{itemize}
\item $e_i = y_{2i-1}y_{2i}x_{i+1}$ for $i \in [3]$, where $x_4 = x_1$,
\item $f_i = x_{\lceil i/2 \rceil} y_i z_i$ for $i \in [6]$, and
\item $g = x_1 x_2 x_3$.
\end{itemize}
\end{construct}

The key vertices of $G_{CP}$, when playing either the Maker-Breaker game or the Chooser-Picker game on $G_{CP}$, are the vertices of $X$.
As these vertices have high degree and are in a sense `central' to $G_{CP}$, both players will be aiming to claim these vertices in their optimal strategies.
Observe that in the Maker-Breaker game, if Maker claims $x_2$ and Breaker claims $x_1$ (say), then Maker has an easy win: He simply claims first $y_3$ (forcing Breaker to claim $z_3$), then at $y_4$ (forcing $z_4$) and finally at $x_3$, claiming the edge $e_2$.
On the other hand if the positions are reversed (that is, if Maker claims $x_1$ and Breaker $x_2$), then the analogous play, that of taking $y_1$ and $y_2$, will fail because Breaker has already claimed $x_2$.
Indeed, if Maker claims $x_1$ then claiming $x_2$ is Breaker's only winning response; in general, if Maker claims a vertex of $X$ then Breaker must claim the vertex of $X$ which lies anticlockwise from it (according to Figure~\ref{CPfigure}).

This asymmetry plays a key role in both the games under consideration. 
In fact, as noted by Bednarska-Bzd\c{e}ga (personal communication), it turns out to be a necessary condition for the results of the Maker-Breaker and Chooser-Picker games on $G_{CP}$ to differ.
In the Chooser-Picker game, this setup favours Chooser; whenever Picker offers a pair of vertices from $X$, Chooser can always force Picker to accept the vertex which lies clockwise from the other.
However, in the Maker-Breaker game, it acts in Breaker's favour since he can see where Maker plays and choose the vertex which lies anticlockwise from it.

\begin{prop} \label{MBProp}
Breaker wins the Maker-Breaker game on $G_{CP}$, where Maker goes first. 
\end{prop}

\proof Without loss of generality we may suppose that Maker claims a vertex of $f_1 \cup f_2$. Breaker the claims $x_2$; this eliminates the edges $e_1, f_3, f_4$ and $g$. Depending on which vertex Maker initially claimed, our analysis splits into two cases:

\medskip \noindent \textbf{Case 1:} \emph{Maker claimed $x_1$.} By Lemma~\ref{Reduction1} we may disregard $e_2$, $f_1$ and $f_2$. 
This leaves only the edges $f_5$, $f_6$ and $e_3$. 
If Maker now claims one of $y_5$ and $y_6$ then Breaker claims the other; otherwise, he claims $x_3$. 
In either case it is clear that he wins.

\medskip \noindent \textbf{Case 2:} \emph{Maker claimed $y_1, y_2, z_1$ or $z_2$.} Without loss of generality Maker claims $y_1$ or $z_1$.
By Lemma~\ref{Reduction1} we may disregard $e_2$ and $f_2$. 
Further, by Lemma~\ref{Reduction2} we may assume that Maker's next move is at $x_1$. 
Then Breaker responds at $y_1$ or $z_1$ (whichever is unclaimed) and eliminates $f_1$. 
Again only the edges $f_5$, $f_6$ and $e_3$ remain and Breaker wins. \endproof

We now turn our attention to the Chooser-Picker game on $G_{CP}$.
Here Picker's most obvious move is to offer two vertices of $X$; as mentioned earlier, however, this attempt fails (see Case 1 below).
Since Picker has a wide range of possible first moves, analysing the game as a whole requires some case analysis.
The general theme in Cases 2--7 is that Chooser will tend to select the vertex of higher degree (when applicable) and will try to play so as to secure a vertex of $X$ for himself, without allowing Picker to claim one in return.

\begin{prop} \label{CPProp} 
Chooser wins the Chooser-Picker game on $G_{CP}$.
\end{prop}

\proof Chooser's winning strategy is as follows. We split into several cases based on Picker's first move:

\medskip \noindent \textbf{Case 1:} \emph{Picker offers two vertices of $X$.} Without loss of generality Picker offers $x_1$ and $x_2$. Then Chooser selects $x_2$. By applying Lemma~\ref{Reduction3} to $f_3$ and $f_4$ we may assume that Picker offers first $\{y_3, z_3\}$, upon which Chooser selects $y_3$, and then $\{y_4, z_4\}$, upon which Chooser selects $y_4$. Now at some later point Chooser claims $x_3$, completing the edge $e_2$ and winning.

\medskip \noindent \textbf{Case 2:} \emph{Picker offers a vertex $x_i$ and a vertex of $f_{2i-1} \cup f_{2i}$.} Without loss of generality $i = 1$. Then Chooser selects $x_1$. By applying Lemma~\ref{Reduction3} to $e_3$ and $g$ we may assume that Picker next offers first $\{y_5, y_6\}$, upon which Chooser selects $y_5$, and then $\{x_2, x_3\}$, upon which Chooser selects $x_3$. Now at some later point Chooser claims $z_5$, completing the edge $f_5$ and winning.

\medskip \noindent \textbf{Case 3:} \emph{Picker offers a vertex $x_i$ and any vertex not in $X$ or $f_{2i-1} \cup f_{2i}$.} Without loss of generality $i = 2$. Then regardless of what the other vertex was, Chooser selects $x_2$. From here we proceed as in Case 1.

\medskip \noindent \textbf{Case 4:} \emph{Picker offers two vertices in $Y$, which are contained in the same edge.} Without loss of generality Picker offers $y_1$ and $y_2$. Then Chooser selects $y_1$. By applying Lemma~\ref{Reduction3} to $f_1$ we may assume that Picker next offers $\{x_1, z_1\}$, upon which Chooser selects $x_1$. From here we proceed as in Case 2.

\medskip \noindent \textbf{Case 5:} \emph{Picker offers a vertex in $Y$ and a vertex in $Z$ which are contained in the same edge.} Without loss of generality Picker offers $y_1$ and $z_1$. Then Chooser selects $y_1$. By applying Lemma~\ref{Reduction3} to $e_1$ we may assume that Picker next offers $\{x_2, y_2\}$, upon which Chooser selects $x_2$. From here we proceed as in Case 1.

\medskip \noindent \textbf{Case 6:} \emph{Picker offers a vertex in $Y$ and a vertex in $Y \cup Z$, which are not contained in the same edge.} Without loss of generality the first vertex Picker offers is $y_1$. Regardless of what the other vertex is, Chooser selects $y_1$. By applying Lemma~\ref{Reduction3} to $f_1$ and $e_1$ we may assume that Picker next offers first $\{x_1, z_1\}$, upon which Chooser selects $x_1$, and then $\{x_2, y_2\}$, upon which Chooser selects $x_2$. Now at some later point Chooser claims $x_3$, claiming the edge $g$ and winning.

\medskip \noindent \textbf{Case 7:} \emph{Picker offers two vertices of $Z$.} Without loss of generality Picker offers either $\{z_1, z_2\}$ or $\{z_1, z_3\}$. In either case Chooser selects $z_1$. By applying Lemma~\ref{Reduction3} to $f_1$ we may assume that Picker next offers $\{x_1, y_1\}$, upon which Chooser selects $x_1$. From here we proceed as in Case 2. \endproof

\section{Construction of $G_4$} \label{G4Section}

The aim of this section is to construct the $4$-graph $G_4$ and to show that Maker wins the Maker-Breaker game on $G_4$. 
Since $G_4$ is rather large and intricate, our construction proceeds in three stages.
We first give a $3$-graph $\Gamma$ which contains the main essential structure which allows Maker to win. 
We then derive a further hypergraph $\Gamma'$ from $\Gamma$ and finally derive $G_4$ from $\Gamma'$. 
Each of these transformations is relatively straightforward; in effect, we successively replace $3$-edges of $\Gamma$ by collections of $4$-edges (which in general involves adding extra vertices to $\Gamma$), while preserving the low maximum degree and the property that Maker wins the Maker-Breaker game on each graph.

\subsection{Construction of $\Gamma$}

\begin{figure}
\begin{center}
\begin{tikzpicture}

\node[label=180:$w_1$] (x11) at (0: 2) {};
\node[label=below:$w_2$] (x21) at (72: 2) {};
\node[label=324:$w_{3}$] (x31) at (144: 2) {};
\node[label=36:$w_{4}$] (x41) at (216: 2) {};
\node[label=above:$w_{5}$] (x51) at (288: 2) {};

\node[label=342:$x_{13}$] (x12) at (342: 5) {};
\node[label=18:$x_{12}$] (x13) at (18: 5) {};
\node[label=54:$x_{23}$] (x22) at (54: 5) {};
\node[label=90:$x_{22}$] (x23) at (90: 5) {};
\node[label=126:$x_{33}$] (x32) at (126: 5) {};
\node[label=162:$x_{32}$] (x33) at (162: 5) {};
\node[label=198:$x_{43}$] (x42) at (198: 5) {};
\node[label=234:$x_{42}$] (x43) at (234: 5) {};
\node[label=270:$x_{53}$] (x52) at (270: 5) {};
\node[label=306:$x_{52}$] (x53) at (306: 5) {};

\node[label=180:$x_{11}$] (x14) at (barycentric cs:x33=0.5,x42=0.5) {};
\node[label=252:$x_{21}$] (x24) at (barycentric cs:x43=0.5,x52=0.5) {};
\node[label=324:$x_{31}$] (x34) at (barycentric cs:x53=0.5,x12=0.5) {};
\node[label=36:$x_{41}$] (x44) at (barycentric cs:x13=0.5,x22=0.5) {};
\node[label=108:$x_{51}$] (x54) at (barycentric cs:x23=0.5,x32=0.5) {};

\node at (189:5.5) {\textcolor{red}{$\mathbf{e_1}$}};
\node at (261:5.5) {\textcolor{red}{$\mathbf{e_2}$}};
\node at (333:5.5) {\textcolor{red}{$\mathbf{e_3}$}};
\node at (45:5.5) {\textcolor{red}{$\mathbf{e_4}$}};
\node at (117:5.5) {\textcolor{red}{$\mathbf{e_5}$}};

\node (t12) at (barycentric cs:x11=0.5,x12=0.5) {};
\node (t22) at (barycentric cs:x21=0.5,x22=0.5) {};
\node (t32) at (barycentric cs:x31=0.5,x32=0.5) {};
\node (t42) at (barycentric cs:x41=0.5,x42=0.5) {};
\node (t52) at (barycentric cs:x51=0.5,x52=0.5) {};

\node (t13) at (barycentric cs:x11=0.5,x13=0.5) {};
\node (t23) at (barycentric cs:x21=0.5,x23=0.5) {};
\node (t33) at (barycentric cs:x31=0.5,x33=0.5) {};
\node (t43) at (barycentric cs:x41=0.5,x43=0.5) {};
\node (t53) at (barycentric cs:x51=0.5,x53=0.5) {};

\node (t14) at (barycentric cs:x11=0.25,x14=0.75) {};
\node (t24) at (barycentric cs:x21=0.25,x24=0.75) {};
\node (t34) at (barycentric cs:x31=0.25,x34=0.75) {};
\node (t44) at (barycentric cs:x41=0.25,x44=0.75) {};
\node (t54) at (barycentric cs:x51=0.25,x54=0.75) {};

%\draw[cap=round, double distance=10] (162:5) -- (198:5);
%\draw[cap=round, double distance=10] (234:5) -- (270:5);
%\draw[cap=round, double distance=10] (306:5) -- (342:5);
%\draw[cap=round, double distance=10] (18:5) -- (54:5);
%\draw[cap=round, double distance=10] (90:5) -- (126:5);

%\draw[cap=round, double distance=10] (162:5) -- (144:2);
%\draw[cap=round, double distance=10] (234:5) -- (216:2);
%\draw[cap=round, double distance=10] (306:5) -- (288:2);
%\draw[cap=round, double distance=10] (18:5) -- (0:2);
%\draw[cap=round, double distance=10] (90:5) -- (72:2);

%\draw[cap=round, double distance=10] (216:2) -- (198:5);
%\draw[cap=round, double distance=10] (288:2) -- (270:5);
%\draw[cap=round, double distance=10] (0:2) -- (342:5);
%\draw[cap=round, double distance=10] (72:2) -- (54:5);
%\draw[cap=round, double distance=10] (144:2) -- (126:5);

\node[draw,color=red,rounded corners=8,fit=(x33) (x14) (x42)] {};
\node[draw,color=red,rounded corners=8,fit=(x33) (x14) (x42),transform canvas={rotate around={72:(0,0)}}] {};
\node[draw,color=red,rounded corners=8,fit=(x33) (x14) (x42),transform canvas={rotate around={144:(0,0)}}] {};
\node[draw,color=red,rounded corners=8,fit=(x33) (x14) (x42),transform canvas={rotate around={216:(0,0)}}] {};
\node[draw,color=red,rounded corners=8,fit=(x33) (x14) (x42),transform canvas={rotate around={288:(0,0)}}] {};

\node[draw,rounded corners=8,fit=(x11) (t14) (x14)] {};
\node[draw,rounded corners=8,fit=(x11) (t14) (x14),transform canvas={rotate around={72:(0,0)}}] {};
\node[draw,rounded corners=8,fit=(x11) (t14) (x14),transform canvas={rotate around={144:(0,0)}}] {};
\node[draw,rounded corners=8,fit=(x11) (t14) (x14),transform canvas={rotate around={216:(0,0)}}] {};
\node[draw,rounded corners=8,fit=(x11) (t14) (x14),transform canvas={rotate around={288:(0,0)}}] {};

\node (aux1) at (11.2825:5) {};
\node (aux2) at (29.2825:2) {};
\node[draw,rounded corners=8,fit=(aux1) (aux2),transform canvas={rotate around={21.3687:(0,0)}}] {};
\node[draw,rounded corners=8,fit=(aux1) (aux2),transform canvas={rotate around={57.3687:(0,0)}}] {};
\node[draw,rounded corners=8,fit=(aux1) (aux2),transform canvas={rotate around={93.3687:(0,0)}}] {};
\node[draw,rounded corners=8,fit=(aux1) (aux2),transform canvas={rotate around={129.3687:(0,0)}}] {};
\node[draw,rounded corners=8,fit=(aux1) (aux2),transform canvas={rotate around={165.3687:(0,0)}}] {};

\node (aux3) at (-11.2825:5) {};
\node (aux4) at (-29.2825:2) {};
\node[draw,rounded corners=8,fit=(aux3) (aux4),transform canvas={rotate around={-21.3687:(0,0)}}] {};
\node[draw,rounded corners=8,fit=(aux3) (aux4),transform canvas={rotate around={-57.3687:(0,0)}}] {};
\node[draw,rounded corners=8,fit=(aux3) (aux4),transform canvas={rotate around={-93.3687:(0,0)}}] {};
\node[draw,rounded corners=8,fit=(aux3) (aux4),transform canvas={rotate around={-129.3687:(0,0)}}] {};
\node[draw,rounded corners=8,fit=(aux3) (aux4),transform canvas={rotate around={-165.3687:(0,0)}}] {};

%\draw[cap=round, double distance=10] (0:2) -- (barycentric cs:x33=0.5,x42=0.5);
%\draw[cap=round, double distance=10] (72:2) -- (barycentric cs:x43=0.5,x52=0.5);
%\draw[cap=round, double distance=10] (144:2) -- (barycentric cs:x53=0.5,x12=0.5);
%\draw[cap=round, double distance=10] (216:2) -- (barycentric cs:x13=0.5,x22=0.5);
%\draw[cap=round, double distance=10] (288:2) -- (barycentric cs:x23=0.5,x32=0.5);

\foreach \i in {1,2,3,4,5} {
        \fill (x\i1) circle (0.1);
}

\foreach \i in {1,2,3,4,5} {
        \fill (x\i2) circle (0.1);
}

\foreach \i in {1,2,3,4,5} {
        \fill (x\i3) circle (0.1);
}

\foreach \i in {1,2,3,4,5} {
        \fill (x\i4) circle (0.1);
}

\foreach \i in {1,2,3,4,5} {
        \fill (t\i2) circle (0.1);
}
\foreach \i in {1,2,3,4,5} {
        \fill (t\i3) circle (0.1);
}
\foreach \i in {1,2,3,4,5} {
        \fill (t\i4) circle (0.1);
}

\end{tikzpicture} 
\caption{The graph $\Gamma$, with the vertices in $W$ and $X$ and the edges $e_1, \ldots, e_5$ labelled.}\label{wGfigure}
\end{center}
\end{figure}
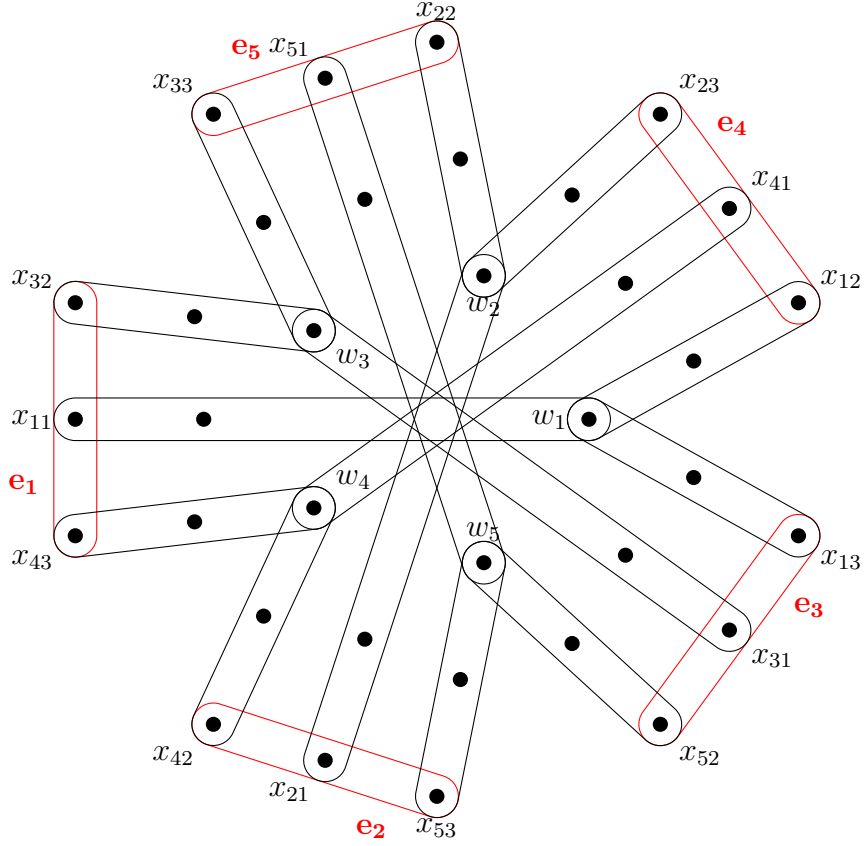

The $3$-graph $\Gamma$ is constructed as follows: Let $W = \{w_{i} \mid i~\in~[5]\}$, $X = \{x_{ij} \mid i~\in~[5], j~\in~[3]\}$ and $T = \{t_{ij} \mid i \in [5], j \in [3]\}$. Then
$$V(\Gamma) = W \cup X \cup T.$$
Let $e_i = x_{i1}x_{(i+2)2}x_{(i+3)3}$ for each $i \in [5]$. 
Then the edges of $\Gamma$ are 
$$\{w_{i}x_{ij}t_{ij} \mid i~\in~[5], j~\in~[3]\} \cup \{e_1, \ldots, e_5\}.$$
(See Figure~\ref{wGfigure}.) 
In total $\Gamma$ has $35$ vertices and $20$ edges. Note that the vertices $w_i$ have degree $3$ for each $i$, the vertices $t_{ij}$ have degree $1$ and the remaining vertices have degree $2$.

%Define the game $Aux$ on $\Gamma$ as follows: $Aux$ has two players, Maker and Breaker (though strictly speaking, $Aux$ is not a Maker-Breaker game in the usual sense), who take turns to claim vertices with Breaker going first. If Maker claims an edge of size $3$ he wins; if he claims an edge of size $2$, he does not immediately win but may immediately claim an additional vertex. Note that we could turn $Aux$ into a genuine Maker-Breaker game by adding an extra vertex to each edge of size $2$ in $\Gamma$, which is not contained in any other edges.

\begin{lemma} \label{weirdGame} Maker has a winning strategy for the Maker-Breaker game on $\Gamma$, where Breaker goes first.
\end{lemma}

\proof Roughly speaking, Maker's strategy is to gain an advantage by playing on the vertices $\{w_i \mid i \in [5]\}$ (the `inner' vertices) and then to use this advantage to claim one of the edges $e_1, \ldots, e_5$. 
By claiming a vertex $w_i$ when $x_{ij}$ and $t_{ij}$ are unclaimed for each $j \in [3]$, Maker effectively gets two `outer' vertices in a single turn: Breaker can claim one of the $x_{ij}$, but Maker will then claim the other two and force a response at the corresponding $t_{ij}$ each time. 
The effect is even greater if Maker has already claimed one of the $x_{ij}$; in this case, Breaker is forced to respond immediately at $t_{ij}$ and Maker gets yet more vertices effectively for free. 
By contrast, if Breaker claims an inner vertex $w_i$ he accomplishes comparatively little. 
Maker was not planning to win on the edges incident to $w_i$ anyway, and he can simply use the above tactic with a different inner vertex.

We now give a rigorous description of Maker's winning strategy for the Maker-Breaker game on $\Gamma$, where Breaker goes first. Whenever Breaker claims a vertex $t_{ij}$ and Maker has not claimed both $w_i$ and $x_{ij}$, Maker plays as if Breaker had claimed $w_i$ or $x_{ij}$, whichever is free (or arbitrarily if neither is free). In other respects the winning strategy is as follows: \medskip

\medskip \noindent \textbf{Case 1:} Breaker claims $w_i$ for some $i$. Without loss of generality Breaker claims $w_1$. 
Maker then claims $w_2$. 
The strategy then branches depending on Breaker's next move:

\medskip \noindent \textbf{Case 1.1:} Breaker claims a vertex of $e_4$ or at $w_4$. 
Then Maker claims $x_{21}$ (forcing $t_{21}$), then $x_{22}$ (forcing $t_{22}$) and then $x_{51}$. 
Breaker is forced to claim $x_{33}$ or Maker wins immediately by playing there and claiming $e_5$. 
Then Maker claims $w_5$ (forcing $t_{51}$), $x_{53}$ (forcing $t_{53}$), and then $x_{42}$, claiming the edge $e_2$ and winning.

\medskip \noindent \textbf{Case 1.2:} Breaker claims a vertex of $e_2$. 
Then Maker claims $x_{22}$ (forcing $t_{22}$), $x_{23}$ (forcing $t_{23}$) and then $x_{41}$. 
Breaker is forced to claim $x_{12}$ or Maker wins immediately by playing there and claiming $e_4$. 
Then Maker claims $w_4$ (forcing $t_{41}$), $x_{43}$ (forcing $t_{43}$) and then  $x_{32}$. 
Breaker is forced to claim $x_{11}$ or Maker wins immediately by playing there and claiming $e_1$. 
Then Maker claims $w_3$ (forcing $t_{32}$), $x_{33}$ (forcing $t_{33}$), and then $x_{51}$, claiming the edge $e_5$ and winning.

\medskip \noindent \textbf{Case 1.3:} Breaker claims any other unclaimed vertex of $\Gamma$. 
Then Maker claims $x_{21}$ (forcing Breaker to claim $t_{21}$), then at $x_{23}$ (forcing $t_{23}$) and then $x_{41}$. 
Breaker is forced to claim $x_{12}$ or Maker wins immediately by playing there and claiming $e_4$. 
Then Maker claims $w_4$ (forcing $t_{41}$), $x_{42}$ (forcing $t_{42}$), and then $x_{53}$, claiming the edge $e_2$ and winning.

\medskip \noindent \textbf{Case 2:} Breaker claims $x_{i1}$ for some $i$. 
Without loss of generality Breaker claims $x_{11}$. 
Maker then claims $w_2$. 
The strategy then branches depending on Breaker's next move:

\medskip \noindent \textbf{Case 2.1:} Breaker claims a vertex of $e_4$ or $w_4$. 
Then Maker wins as in Case 1.1, which is possible since Breaker's first move at $x_{11}$ does not prevent this.

\medskip \noindent \textbf{Case 2.2:} Breaker claims a vertex of $e_2$. 
Then Maker claims $x_{22}$ (forcing $t_{22}$), $x_{23}$ (forcing $t_{23}$) and then $x_{12}$. 
Breaker is forced to claim $x_{41}$ or Maker wins immediately by playing there and claiming $e_4$. 
Then Maker claims $w_1$ (forcing $t_{12}$), $x_{13}$ (forcing $t_{13}$) and then $x_{52}$. 
Breaker is forced to claim $x_{31}$ or Maker wins immediately by playing there and claiming $e_3$. 
Then Maker claims $w_5$ (forcing $t_{52}$), $x_{51}$ (forcing $t_{51}$) and then $x_{33}$, claiming the edge $e_5$ and winning.

\medskip \noindent \textbf{Case 2.3:} Breaker claims any other unclaimed vertex. Then Maker wins as in Case 1.3.

\medskip \noindent \textbf{Case 3:} Breaker claims $x_{i2}$ or $x_{i3}$ for some $i$. Without loss of generality Breaker claims $x_{32}$. Maker then claims $w_2$. The strategy then branches depending on Breaker's next move:

\medskip \noindent \textbf{Case 3.1:} Breaker claims a vertex of $e_4$ or $w_4$. Then Maker wins as in Case 1.1.

\medskip \noindent \textbf{Case 3.2:} Breaker claims a vertex of $e_2$. Then Maker wins as in Case 2.2. 

\medskip \noindent \textbf{Case 3.3:} Breaker claims any other unclaimed vertex. Then Maker wins as in Case 1.3. \endproof

\subsection{Construction of $\Gamma'$ from $\Gamma$}

\begin{figure}
\begin{center}
\begin{tikzpicture}[label distance=2]

\node[label=left:${w=w_i}$] (ox11) at (-3,2) {};
\fill (ox11) circle (0.1);
\node[label=left:${x=x_{ij}}$] (ox12) at (-3,-2) {};
\fill (ox12) circle (0.1);
\node[label=left:${t=t_{ij}}$] (ot12) at (-3,0) {};
\fill (ot12) circle (0.1);

\node[label=below left:$w$] (x11) at (0, 2) {};
\fill (x11) circle (0.1);
\node[label=above left:$x$] (x12) at (0, -2) {};
\fill (x12) circle (0.1);
\node[label=above left:$t$] (t12) at (2, 0) {};
\fill (t12) circle (0.1);

\begin{scope}[label distance=2]
\node[label=right:$y_{1}$] (y121) at (3,-1) {};
\node[label=right:$y_{2}$] (y122) at (4.5,-1) {};
\node[label=right:$y_{3}$] (y123) at (3,0) {};
\node[label=right:$y_{4}$] (y124) at (4.5,0) {};
\node[label=right:$y_{5}$] (y125) at (3,1) {};
\node[label=right:$y_{6}$] (y126) at (4.5,1) {};
\end{scope}

\foreach \i in {1,2,3,4,5,6} {
        \fill (y12\i) circle (0.1);
}

\node[label=below:$z_{1}$] (z121) at (2,-2) {};
\node[label=below:$z_{2}$] (z122) at (3,-2) {};
\node[label=below:$z_{3}$] (z123) at (4.5,-2) {};
\node[label=below:$z_{4}$] (z124) at (5.5,-2) {};

\foreach \i in {1,2,3,4} {
        \fill (z12\i) circle (0.1);
}

\draw[rounded corners=8] ($(ox11) + (-0.2828, 0.2828)$) -- ($(ox11) + (0.2828, 0.2828)$) -- ($(ox12) + (0.2828, -0.2828)$) -- ($(ox12) + (-0.2828, -0.2828)$) -- cycle;
\draw[->, line width=5] (-2, 0) -- (-0.5, 0);

\draw[rounded corners=8] ($(x11) + (-0.4,0)$) -- ($(y121) + (-0.1112,-0.2828)$) -- ($(y122) + (0.2828,-0.2828)$) -- ($(y122) + (0.2828,0.2828)$) -- ($(y121) + (0.1112,0.2828)$) -- ($(x11) + (0,0.4)$) -- cycle;
\draw[rounded corners=8] ($(x12) + (-0.4,0)$) -- ($(y125) + (-0.1112,0.2828)$) -- ($(y126) + (0.2828,0.2828)$) -- ($(y126) + (0.2828,-0.2828)$) -- ($(y125) + (0.1112,-0.2828)$) -- ($(x12) + (0,-0.4)$) -- cycle;
\draw[rounded corners=8] ($(x12) + (-0.4,0)$) -- ($(t12) + (-0.1112,0.2828)$) -- ($(y124) + (0.2828,0.2828)$) -- ($(y124) + (0.2828,-0.2828)$) -- ($(t12) + (0.1112,-0.2828)$) -- ($(x12) + (0,-0.4)$) -- cycle;
\draw[rounded corners=8] ($(z121) + (-0.4,0)$) -- ($(y121) + (-0.2828,0.1112)$) -- ($(y125) + (-0.2828,0.2828)$) -- ($(y125) + (0.2828,0.2828)$) -- ($(y121) + (0.2828,-0.1112)$) -- ($(z121) + (0,-0.4)$) -- cycle;
\draw[rounded corners=8] ($(z122) + (-0.2828,-0.2828)$) -- ($(y121) + (-0.2828,-0.1112)$) -- ($(y125) + (-0.2828,0.2828)$) -- ($(y125) + (0.2828,0.2828)$) -- ($(y121) + (0.2828,-0.1112)$) -- ($(z122) + (0.2828,-0.2828)$) -- cycle;
\draw[rounded corners=8] ($(z123) + (-0.2828,-0.2828)$) -- ($(y126) + (-0.2828,0.2828)$) -- ($(y126) + (0.2828,0.2828)$) -- ($(z123) + (0.2828,-0.2828)$) -- cycle;
\draw[rounded corners=8] ($(z124) + (0.4,0)$) -- ($(y122) + (0.2828,0.1112)$) -- ($(y126) + (0.2828,0.2828)$) -- ($(y126) + (-0.2828,0.2828)$) -- ($(y122) + (-0.2828,-0.1112)$) -- ($(z124) + (0,-0.4)$) -- cycle;

\end{tikzpicture}
\caption{Forming $\Gamma'$ from $\Gamma$.}\label{nGfigure}
\end{center}
\end{figure}
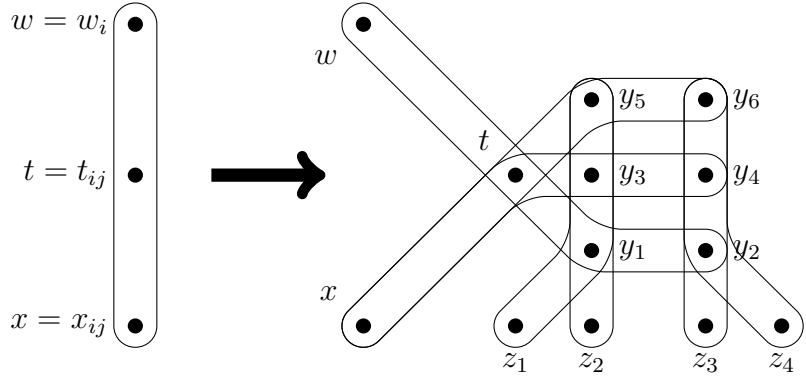

We derive the (non-uniform) hypergraph $\Gamma'$ from $\Gamma$ as follows: Add new vertices $(y_{ijk})_{k \in [6]}$ and $(z_{ijk})_{k \in [4]}$ to $\Gamma$ for every $i \in [5]$ and $j \in [3]$. 
Then for every $i \in [5]$ and $j \in [3]$, writing $w = w_i$, $x = x_{ij}$, $t = t_{ij}$, $y_k = y_{ijk}$ for $k \in [6]$ and $z_k = z_{ijk}$ for $k \in [4]$, we replace the edge $wxt$ by the edges 
$$w ty_1y_{2},\, xty_{3}y_{4},\, w ty_{5}y_{6}$$
and
$$y_{1}y_{3}y_{5}z_{1},\, y_{1}y_{3}y_{5}z_{2},\, y_{2}y_{4}y_{6}z_{3},\, y_{2}y_{4}y_{6}z_{4}.$$ (See Figure~\ref{nGfigure}.)

In total $\Gamma'$ has $155$ vertices and $110$ edges. Note that all of the edges of $\Gamma'$, apart from the edges $e_1, \ldots, e_5$, have size $4$, since all of the other original edges of $\Gamma$ have been replaced and every edge we added had size $4$. Further, $\Gamma'$ has maximum degree $3$. Indeed, at each of the vertices of $W$ we replaced each $3$-edge by a single $4$-edge, at the vertices of $X$ we left one $3$-edge as it is and replaced the other by two $4$-edges, and Figure~\ref{nGfigure} illustrates that the vertices $z_{ijk}$ have degree $1$ and the remaining vertices degree $3$.

\begin{lemma} \label{normalGame} Maker has a winning strategy for the Maker-Breaker game on $\Gamma'$, where Breaker goes first.
\end{lemma}

\proof Initially Maker plays only on the vertices $\{w_i \mid i \in [5]\}$ and $\{x_{ij} \mid i \in [5], j \in [3]\}$ and plays according to the winning strategy for $\Gamma$ from Lemma \ref{weirdGame}. 
If Breaker claims $y_{ijk}$ or $z_{ijk}$ for some $i, j, k$ then Maker plays as if Breaker had claimed $t_{ij}$.
(If Breaker had previously claimed $t_{ij}$ then we choose a free $x_{ij}$ arbitrarily and play as if Breaker had played there instead.)

Since Maker is following the winning strategy for $\Gamma$, at some point he will claim an edge of $\Gamma$. 
If this edge is also an edge of $\Gamma'$ then he wins immediately, so we may assume that Maker claims $w_i x_{ij}t_{ij}$ for some $i \in [5]$ and $j \in [3]$. 
Since Maker claimed $t_{ij}$, Breaker cannot have claimed $y_{ijk}$ or $z_{ijk}$ for any $k$. 
Since moves by Breaker apart from $y_{ijk}$ or $z_{ijk}$ for some $k$ will not affect what follows, we may assume that Breaker's next move is at one of these vertices.
Suppose that Breaker's next move is either $y_{ij1}$, $z_{ij1}$ or $z_{ij2}$ (other possible moves are covered in a similar way).
Now Maker claims $y_{ij4}$. Breaker is forced to claim $y_{ij3}$, or Maker plays there and wins. 
Maker then claims $y_{ij6}$, forcing Breaker to claim $y_{ij5}$. Finally Maker claims $y_{ij2}$, and now takes one of $z_{ij3}$ and $z_{ij4}$ and wins. \endproof

\subsection{Construction of $G_4$ from $\Gamma'$}

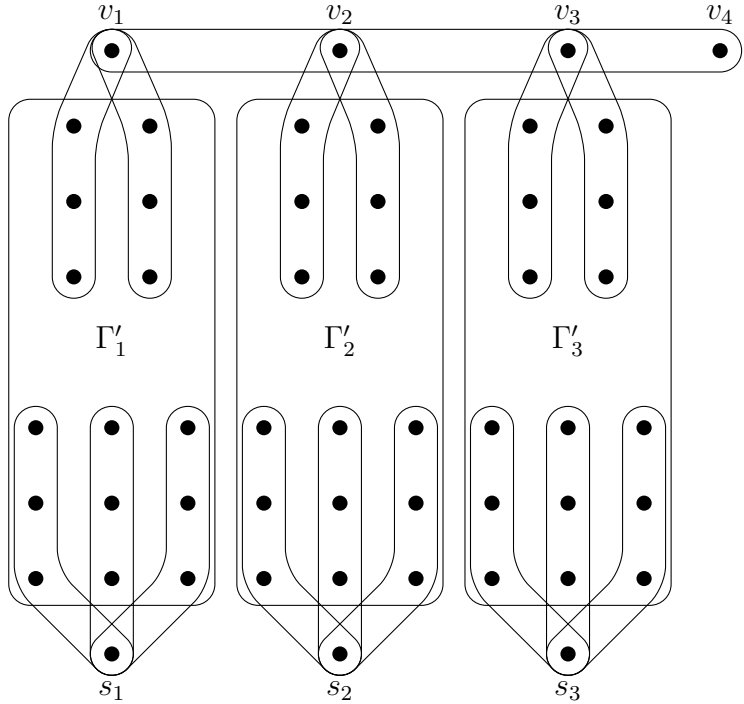
\begin{figure}
\begin{center}
\begin{tikzpicture}[label distance=2]

\node[label=above:$v_1$] (v1) at (2, 9) {};
\node[label=above:$v_2$] (v2) at (5, 9) {};
\node[label=above:$v_3$] (v3) at (8, 9) {};
\node[label=above:$v_4$] (v4) at (10, 9) {};

\foreach \i in {1,2,3,4} {
        \fill (v\i) circle (0.1);
}

\node[label=below:$s_1$] (s1) at (2, 1) {};
\node[label=below:$s_2$] (s2) at (5, 1) {};
\node[label=below:$s_3$] (s3) at (8, 1) {};

\foreach \i in {1,2,3} {
        \fill (s\i) circle (0.1);
}

\node (w11) at (1.5, 8) {};
\node (w12) at (1.5, 7) {};
\node (w13) at (1.5, 6) {};

\node (w21) at (2.5, 8) {};
\node (w22) at (2.5, 7) {};
\node (w23) at (2.5, 6) {};

\node (w31) at (4.5, 8) {};
\node (w32) at (4.5, 7) {};
\node (w33) at (4.5, 6) {};

\node (w41) at (5.5, 8) {};
\node (w42) at (5.5, 7) {};
\node (w43) at (5.5, 6) {};

\node (w51) at (7.5, 8) {};
\node (w52) at (7.5, 7) {};
\node (w53) at (7.5, 6) {};

\node (w61) at (8.5, 8) {};
\node (w62) at (8.5, 7) {};
\node (w63) at (8.5, 6) {};

\node (w71) at (1, 2) {};
\node (w72) at (1, 3) {};
\node (w73) at (1, 4) {};

\node (w81) at (2, 2) {};
\node (w82) at (2, 3) {};
\node (w83) at (2, 4) {};

\node (w91) at (3, 2) {};
\node (w92) at (3, 3) {};
\node (w93) at (3, 4) {};

\foreach \i in {1,2,...,9} {
        \fill (w\i1) circle (0.1);
}

\foreach \i in {1,2,...,9} {
        \fill (w\i2) circle (0.1);
}

\foreach \i in {1,2,...,9} {
        \fill (w\i3) circle (0.1);
}

\node (wa1) at (4, 2) {};
\node (wa2) at (4, 3) {};
\node (wa3) at (4, 4) {};

\foreach \i in {1,2,3} {
        \fill (wa\i) circle (0.1);
}

\node (wb1) at (5, 2) {};
\node (wb2) at (5, 3) {};
\node (wb3) at (5, 4) {};

\foreach \i in {1,2,3} {
        \fill (wb\i) circle (0.1);
}

\node (wc1) at (6, 2) {};
\node (wc2) at (6, 3) {};
\node (wc3) at (6, 4) {};

\foreach \i in {1,2,3} {
        \fill (wc\i) circle (0.1);
}

\node (wd1) at (7, 2) {};
\node (wd2) at (7, 3) {};
\node (wd3) at (7, 4) {};

\foreach \i in {1,2,3} {
        \fill (wd\i) circle (0.1);
}

\node (we1) at (8, 2) {};
\node (we2) at (8, 3) {};
\node (we3) at (8, 4) {};

\foreach \i in {1,2,3} {
        \fill (we\i) circle (0.1);
}

\node (wf1) at (9, 2) {};
\node (wf2) at (9, 3) {};
\node (wf3) at (9, 4) {};

\foreach \i in {1,2,3} {
        \fill (wf\i) circle (0.1);
}

\draw[rounded corners=8] ($(v1) + (-0.2828, 0.2828)$) -- ($(v1) + (-0.2828, -0.2828)$) -- ($(v4) + (0.2828, -0.2828)$) -- ($(v4) + (0.2828, 0.2828)$) -- cycle;

\draw[rounded corners=8] ($(s1) + (-0.4,0)$) -- ($(w91) + (-0.2828,0.1112)$) -- ($(w93) + (-0.2828,0.2828)$) -- ($(w93) + (0.2828,0.2828)$) -- ($(w91) + (0.2828,-0.1112)$) -- ($(s1) + (0,-0.4)$) -- cycle;
\draw[rounded corners=8] ($(s2) + (-0.4,0)$) -- ($(wc1) + (-0.2828,0.1112)$) -- ($(wc3) + (-0.2828,0.2828)$) -- ($(wc3) + (0.2828,0.2828)$) -- ($(wc1) + (0.2828,-0.1112)$) -- ($(s2) + (0,-0.4)$) -- cycle;
\draw[rounded corners=8] ($(s3) + (-0.4,0)$) -- ($(wf1) + (-0.2828,0.1112)$) -- ($(wf3) + (-0.2828,0.2828)$) -- ($(wf3) + (0.2828,0.2828)$) -- ($(wf1) + (0.2828,-0.1112)$) -- ($(s3) + (0,-0.4)$) -- cycle;

\draw[rounded corners=8] ($(s1) + (0.4,0)$) -- ($(w71) + (0.2828,0.1112)$) -- ($(w73) + (0.2828,0.2828)$) -- ($(w73) + (-0.2828,0.2828)$) -- ($(w71) + (-0.2828,-0.1112)$) -- ($(s1) + (0,-0.4)$) -- cycle;
\draw[rounded corners=8] ($(s2) + (0.4,0)$) -- ($(wa1) + (0.2828,0.1112)$) -- ($(wa3) + (0.2828,0.2828)$) -- ($(wa3) + (-0.2828,0.2828)$) -- ($(wa1) + (-0.2828,-0.1112)$) -- ($(s2) + (0,-0.4)$) -- cycle;
\draw[rounded corners=8] ($(s3) + (0.4,0)$) -- ($(wd1) + (0.2828,0.1112)$) -- ($(wd3) + (0.2828,0.2828)$) -- ($(wd3) + (-0.2828,0.2828)$) -- ($(wd1) + (-0.2828,-0.1112)$) -- ($(s3) + (0,-0.4)$) -- cycle;

\draw[rounded corners=8] ($(w83) + (-0.2828, 0.2828)$) -- ($(w83) + (0.2828, 0.2828)$) -- ($(s1) + (0.2828, -0.2828)$) -- ($(s1) + (-0.2828, -0.2828)$) -- cycle;
\draw[rounded corners=8] ($(wb3) + (-0.2828, 0.2828)$) -- ($(wb3) + (0.2828, 0.2828)$) -- ($(s2) + (0.2828, -0.2828)$) -- ($(s2) + (-0.2828, -0.2828)$) -- cycle;
\draw[rounded corners=8] ($(we3) + (-0.2828, 0.2828)$) -- ($(we3) + (0.2828, 0.2828)$) -- ($(s3) + (0.2828, -0.2828)$) -- ($(s3) + (-0.2828, -0.2828)$) -- cycle;

\draw[rounded corners=8] ($(v1) + (-0.3464,0.2)$) -- ($(w21) + (-0.2828,-0.1682)$) -- ($(w23) + (-0.2828,-0.2828)$) -- ($(w23) + (0.2828,-0.2828)$) -- ($(w21) + (0.2828,0)$) -- ($(v1) + (0.2,0.3464)$) -- cycle;
\draw[rounded corners=8] ($(v2) + (-0.3464,0.2)$) -- ($(w41) + (-0.2828,-0.1682)$) -- ($(w43) + (-0.2828,-0.2828)$) -- ($(w43) + (0.2828,-0.2828)$) -- ($(w41) + (0.2828,0)$) -- ($(v2) + (0.2,0.3464)$) -- cycle;
\draw[rounded corners=8] ($(v3) + (-0.3464,0.2)$) -- ($(w61) + (-0.2828,-0.1682)$) -- ($(w63) + (-0.2828,-0.2828)$) -- ($(w63) + (0.2828,-0.2828)$) -- ($(w61) + (0.2828,0)$) -- ($(v3) + (0.2,0.3464)$) -- cycle;

\draw[rounded corners=8] ($(v1) + (0.3464,0.2)$) -- ($(w11) + (0.2828,-0.1682)$) -- ($(w13) + (0.2828,-0.2828)$) -- ($(w13) + (-0.2828,-0.2828)$) -- ($(w11) + (-0.2828,0)$) -- ($(v1) + (-0.2,0.3464)$) -- cycle;
\draw[rounded corners=8] ($(v2) + (0.3464,0.2)$) -- ($(w31) + (0.2828,-0.1682)$) -- ($(w33) + (0.2828,-0.2828)$) -- ($(w33) + (-0.2828,-0.2828)$) -- ($(w31) + (-0.2828,0)$) -- ($(v2) + (-0.2,0.3464)$) -- cycle;
\draw[rounded corners=8] ($(v3) + (0.3464,0.2)$) -- ($(w51) + (0.2828,-0.1682)$) -- ($(w53) + (0.2828,-0.2828)$) -- ($(w53) + (-0.2828,-0.2828)$) -- ($(w51) + (-0.2828,0)$) -- ($(v3) + (-0.2,0.3464)$) -- cycle;

\node[inner sep=6,draw,rounded corners=8,fit=(w11) (w91) (w71)] {$\Gamma'_1$};
\node[inner sep=6,draw,rounded corners=8,fit=(w31) (wc1) (wa1)] {$\Gamma'_2$};
\node[inner sep=6,draw,rounded corners=8,fit=(w51) (wf1) (wd1)] {$\Gamma'_3$};

\end{tikzpicture}
\caption{Construction of $G_4$, using three copies of $\Gamma'$.} \label{Gfigure}
\end{center}
\end{figure}

Now we construct a $4$-graph $G_4$ as follows: Form the disjoint union of three copies $\Gamma'_1, \Gamma'_2, \Gamma'_3$ of $\Gamma'$, and let $e_{i1}, \ldots, e_{i5}$ be the $3$-edges of $\Gamma'_i$ for each $i \in [3]$. 
Add extra vertices $v_1, v_2, v_3, v_4, s_1, s_2, s_3$. Now replace the edges $e_{i1}$ and $e_{i2}$ by $e_{i1} \cup v_i$ and $e_{i2} \cup v_i$, and replace the edges $e_{i3}, e_{i4}, e_{i5}$ by $e_{i3} \cup s_i, e_{i4} \cup s_i, e_{i5} \cup s_i$ for each $i \in [3]$.
Finally add the edge $v_1 v_2 v_3 v_4$.
(See Figure~\ref{Gfigure}.)

Note that $G_4$ is indeed a $4$-graph, since we replaced each $3$-edge $e_{ij}$ by either the $4$-edge $e_{ij} \cup v_i$ or the $4$-edge $e_{ij} \cup s_i$. Further, it is easily seen that $G_4$ has maximum degree $3$. As mentioned previously, in total $G_4$ has $472$ vertices and $331$ edges.

\begin{prop} \label{4GraphGame} Maker has a winning strategy for the Maker-Breaker game on $G_4$, where Maker goes first.
\end{prop}

\proof Firstly, suppose that at some point during the game Maker claims $v_i$ for some $i \in [3]$, such that Breaker has not yet claimed $s_i$ or at any vertex of $\Gamma'_i$. Suppose further that Breaker does not immediately respond at any of these vertices. Then Maker claims $s_i$ and can now use the Breaker-first winning strategy on $\Gamma'_i$ from Lemma \ref{normalGame}. 
Thus Maker claims an edge of $\Gamma'_i$.
But any $4$-edge of $\Gamma'_i$ is also an edge of $G_4$ and if Maker claims a $3$-edge $e_{ij}$ then he also claims either the $4$-edge $e_{ij} \cup v_i$ or the $4$-edge $e_{ij} \cup s_i$ (whichever is an edge of $G_4$).
Hence Maker wins in either case. 

So we may assume that every time Maker claims $v_{i}$ for some $i \in [3]$, Breaker immediately responds at $s_i$ or at some vertex of $\Gamma'_i$, assuming he has not done so already. 
But in this case, Maker simply claims each $v_{i}$ in turn, ending with $v_4$, and wins. \endproof

\end{document}